\documentclass[11pt,reqno]{article}
\usepackage{latexsym, amsmath, amssymb, amsthm, a4, epsfig}
\usepackage{graphicx}
\usepackage{color}

\newtheorem{theorem}{Theorem}[section]

\newtheorem{lemma}[theorem]{Lemma}
\newtheorem{prop}[theorem]{Proposition}

\newcommand{\p}{\partial}

\newcommand{\eqnref}[1]{(\ref {#1})}

\newcommand{\Rbb}{\mathbb{R}}

\newcommand{\la}{\langle}
\newcommand{\ra}{\rangle}

\newcommand{\Fcal}{\mathcal{F}}

\newcommand{\Kcal}{\mathcal{K}}

\newcommand{\Gd}{\delta}

\newcommand{\Gvf}{\varphi}

\newcommand{\Gl}{\lambda}

\newcommand{\Gs}{\sigma}

\newcommand{\Gz}{\zeta}

\newcommand{\GG}{\Gamma}

\newcommand{\GO}{\Omega}

\newcommand{\beq}{\begin{equation}}
\newcommand{\eeq}{\end{equation}}

\def\ol{\overline}

\numberwithin{equation}{section}
\numberwithin{figure}{section}

\begin{document}

\title{Spectral structure of the Neumann--Poincar\'e operator on thin domains in two dimensions\thanks{\footnotesize This work was supported by NRF (of S. Korea) grants No. 2019R1A2B5B01069967 and by JSPS (of Japan) KAKENHI Grant Number JP19K14553.}}

\author{Kazunori Ando\thanks{Department of Electrical and Electronic Engineering and Computer Science, Ehime University, Ehime 790-8577, Japan. Email: {\tt ando@cs.ehime-u.ac.jp}.}
\and Hyeonbae Kang\thanks{Department of Mathematics and Institute of Applied Mathematics, Inha University, Incheon 22212, S. Korea. Email: {\tt hbkang@inha.ac.kr}.}
\and Yoshihisa Miyanishi\thanks{Center for Mathematical Modeling and Data Science, Osaka University, Osaka 560-8531, Japan. Email: {\tt miyanishi@sigmath.es.osaka-u.ac.jp}.}}
\date{}
\maketitle

\begin{abstract}
We consider the spectral structure of the Neumann--Poincar\'e operators defined on the boundaries of thin domains of rectangle shape in two dimensions. We prove that as the aspect ratio of the domains tends to $\infty$, or equivalently, as the domains get thinner, the spectra of the Neumann--Poincar\'e operators are densely distributed in the interval $[-1/2,1/2]$.
\end{abstract}

\noindent{\footnotesize {\bf AMS subject classifications}. 35J47 (primary), 35P05 (secondary)}

\noindent{\footnotesize {\bf Key words}. Neumann--Poincar\'e operator, spectrum, thin domain, Poisson kernel}


\section{Introduction and statements of results}

For a bounded domain $\GO$ with the Lipschitz continuous boundary in $\Rbb^2$, the Neumann--Poincar\'e (abbreviated by NP) operator is the boundary integral operator on $\p\GO$ defined by
\beq
\Kcal_{\p\GO} [\Gvf](x) = \mbox{p.v.} \frac{1}{2\pi} \int_{\p\GO_R} \frac{\la y-x, {\bf \nu}_{y} \ra}{|x-y|^2} \Gvf(y) ds(y), \quad x \in \p\GO,
\eeq
where p.v. stands for the Cauchy principal value and $ds$ the line element on $\p\GO$. This operator naturally appears when solving the classical Dirichlet or Neumann problems using layer potentials. The NP operator is also called the double layer potential.

Even though $\Kcal_{\p\GO}$ is not self-adjoint on $L^2(\p\GO)$ unless $\p\GO$ is a disk or a ball, it can be realized as a self-adjoint operator on $H^{1/2}(\p\GO)$ (the $L^2$-Sobolev space of order $1/2$) using the Plemelj's symmetrization principle (see \cite{KPS-ARMA-07}). Thus the spectrum of the NP operator on $H^{1/2}(\p\GO)$, which is denoted by $\Gs(\Kcal_{\p\GO})$, consists of continuous spectrum and pure point spectrum \cite{Yo}.

It is proved lately in \cite{PP2} that if a two-dimensional domain $\GO$ has corners on its boundary, then $\Kcal_{\p\GO}$ has continuous spectrum which is a connected interval symmetric with respect to $0$, and the end points of the interval is completely determined by the smallest angle of the corners. In particular, if $\GO$ is a rectangle, then the continuous spectrum is the interval $[-1/4, 1/4]$. It is known that $\Gs(\Kcal_{\p\GO}) \setminus \{1/2\}$ is the closed subset of $(-1/2,1/2)$. In recent work \cite{HKL}, a classification method to distinguish eigenvalues from continuous spectrum has been proposed, and implemented numerically to investigate existence of eigenvalues on various domains with corners. The numerical experiments reveal that on rectangles more and more eigenvalues of the NP operator appear outside the interval $[-1/4, 1/4]$ of the continuous spectrum as the aspect ratio of the rectangle gets larger. It is also proved that if the aspect ratio is large enough, there is at least one eigenvalue outside $[-1/4, 1/4]$. In this paper we improve this result drastically and prove that the spectra actually fill up the whole interval $(-1/2,1/2)$ in some sense as the aspect ratio gets larger.

The two-dimensional domains to be considered in this paper are not just rectangles. The long sides are lines, but the short sides do not have to be lines, they can be curves. Since the NP operator is dilation invariant, we define planar thin domains as follows: for $R\ge 1$, let $\GO_R$ be either the rectangle of the aspect ratio $R$ or the domain obtained after smoothing the corners of the rectangle: the boundary $\p\GO_R$, which is smooth or Lipschitz continuous allowing corners, consists of three parts, say
\beq\label{boundary}
\p\GO_R = \GG_R^+ \cup \GG_R^- \cup \GG_R^s
\eeq
where the top and bottom are
\beq
\GG_R^+=[-R,R] \times \{1\}, \quad \GG_R^-=[-R,R] \times \{-1\},
\eeq
and the side $\GG_R^s$ consists of the left and right sides, namely, $\GG_R^s=\GG_R^l \cup \GG_R^r$, where $\GG_R^l$ and $\GG_R^r$ are curves connecting points $(\mp R, 1)$ and $(\mp R, -1)$, respectively. We assume that $\GG_R^l$ and $\GG_R^r$ are of any but fixed shape independent of $R$.

The following theorem is the main result of this paper.
\begin{theorem}\label{thm:main}
If $\{ R_j \}$ be an increasing sequence such that $R_j \to \infty$ as $j \to \infty$, then
\beq\label{main}
\ol{\cup_{j=1}^\infty \Gs(\Kcal_{\p\GO_{R_j}})} = [-1/2, 1/2].
\eeq
\end{theorem}

If $\p\GO_R$ is smooth, then eigenvalues of each $\Kcal_{\p\GO_{R_j}}$ accumulates to $0$ since it is a compact operator. However, Theorem \ref{thm:main} shows that more and more eigenvalues are approaching to $\pm 1/2$ as $j \to \infty$, that is, as the domains get longer (or equivalently, thinner). It actually show that the totality of eigenvalues of $\Kcal_{\p\GO_{R_j}}$ are densely distributed in $[-1/2, 1/2]$ as $R_j \to \infty$. It is rather surprising that \eqnref{main} holds regardless of the choice of $\{ R_j \}$. If $\GO_R$ are rectangles, then the totality of eigenvalues of $\Kcal_{\p\GO_{R_j}}$ is dense outside continuous spectrum, namely, in $[-1/2,-1/4] \cup [1/4,1/2]$. In particular, there are infinitely many eigenvalues outside continuous spectrum.

The characteristic feature of the set $\GO_R$ is that it tends to the infinite strip $\Rbb \times [0,2]$ as $R \to \infty$. Furthermore, the NP operator on $\p\GO_R$ behaves like the Poisson integral. Since the Fourier transform of the Poisson kernel is $e^{-2\pi t |\xi|}$, the Poisson integral has $[0,1]$ as its continuous spectrum. This is the key observation in proving Theorem \ref{thm:main}.

This paper is organized as follows. In the next section we relates the NP operator on thin domain with the Poisson integral on the half space. Theorem \ref{thm:main} is proved in section \ref{sec3}. This paper ends with discussions on some related problems.

\section{Poisson integral and construction of test functions}

To motivate the construction of test functions in this section, we define
\beq\label{Gvfdef}
\Gvf(x)= \Gvf(x_1,x_2) :=
\begin{cases}
f(x_1) \quad &\mbox{if } x \in \GG_R^+ \cup \GG_R^-, \\
0 \quad &\mbox{if } x \in \GG_R^s
\end{cases}
\eeq
for a given a compactly supported function $f$ on $[-R,R]$. Let us put $\Kcal_R=\Kcal_{\p\GO_{R}}$ for simplicity of notation, from here on.
Note that
\begin{align}
\Kcal_R[\Gvf](x_1,x_2) & = -\frac{1}{2\pi} \int_{\Rbb^1} \frac{x_2-1}{(x_1-y_1)^2 + (x_2-1)^2} f(y_1) dy_1 \nonumber \\
& \qquad + \frac{1}{2\pi} \int_{\Rbb^1} \frac{x_2+1}{(x_1-y_1)^2 + (x_2+1)^2} f(y_1) dy_1. \label{Poisson}
\end{align}
Thus, if $(x_1,x_2) \in \GG_R^+ \cup \GG_R^-$, namely, if $x_2=1$ or $-1$, then
\beq
\Kcal_R[\Gvf](x_1,x_2) = \frac{1}{2\pi} \int_{\Rbb^1} \frac{2}{(x_1-y_1)^2 + 2^2} f(y_1) dy_1 = \frac{1}{2} P_2 * f(x_1),
\eeq
where $P_t$ is the Poisson kernel on the half space, namely,
\beq
P_t(x_1) = \frac{1}{\pi} \frac{t}{|x_1|^2+t^2}, \quad t>0.
\eeq

We look for a function $f$ such that
\beq\label{GlPf}
\Gl f - \frac{1}{2} P_2 * f =0
\eeq
for a given $\Gl \in (0, 1/2]$. Since $\hat{P_t}(\xi) = \exp(-2\pi t |\xi|)$, it amounts to
\beq\label{zeroeqn}
\left( \Gl - \frac{1}{2} e^{-4\pi |\xi|} \right) \hat{f}(\xi)=0.
\eeq
Here, the Fourier transform is defined by
$$
\hat{f}(\xi)=\Fcal[f](\xi):= \int_{\Rbb^1} e^{-2\pi i \xi x} f(x) dx.
$$
Let $\xi_0 >0$ be such that
\beq\label{GLdef}
\Gl - \frac{1}{2} e^{-4\pi |\xi_0|} =0.
\eeq
Such a point exists since $\Gl \in (0, 1/2]$. The relation \eqnref{zeroeqn} can be satisfied only when $\hat{f}$ is supported on the set $\{ |\xi|=|\xi_0| \}$. Thus $\hat{f}(\xi)= \Gd_{\xi_0}(\xi)$, the Dirac-delta function at $\xi_0$, is a good candidate.

We now construct the desired test functions $f_R$ to be used in the next section.
Let $\psi$ be a function such that $\hat{\psi}$ is a non-negative compactly supported smooth function such that
$$
\int_{\Rbb^1} \hat{\psi}(\xi)d\xi=1.
$$
Then, $R \hat{\psi}(R(\xi-\xi_0))$ converges weakly to $\Gd_{\xi_0}(\xi)$ as $R \to \infty$. Define $g_R$ by
\beq
\widehat{g_R}(\xi)= R \hat{\psi}(R(\xi-\xi_0)).
\eeq
Then one can see easily that
\beq
g_R(x)= e^{2\pi i \xi_0 x} \psi(R^{-1}x).
\eeq
Let $\chi$ be a smooth cut-off function such that $\mbox{supp}(\chi) \subset [-1/2,1/2]$ and $\chi=1$ on $[-1/4, 1/4]$. Define
\beq\label{fRx}
f_R(x):= \chi(R^{-1}x) g_R(x)=  e^{2\pi i \xi_0 x} (\chi\psi)(R^{-1}x).
\eeq

In what follows, $\| \cdot \|_{1/2}$ denotes the Sobolev $1/2$-norm on $\Rbb^1$, and $A \lesssim B$ means that there is a constant independent of $R$ such that $A \le CB$.

\begin{lemma}
For $\Gl \in (0,1/2]$, let $f_R$ be the function defined by \eqnref{fRx} where $\xi_0$ satisfies \eqnref{GLdef}.
It holds that
\beq\label{fR1/2}
R^{1/2} \lesssim \| f_R \|_{1/2}
\eeq
and
\beq\label{Gl-Kest}
\left\| \Gl f_R -\frac{1}{2}P_2 *f_R \right\|_{1/2} \lesssim R^{-1/2}.
\eeq
\end{lemma}

\proof
Note that
\beq
\widehat{f_R}(\xi)= R \widehat{(\chi\psi)}(R(\xi-\xi_0)).
\eeq
Thus,
\begin{align*}
\| f_R \|_{1/2}^2 &= \int_{\Rbb} (1+|\xi|) |\widehat{f_R}(\xi)|^2 d\xi \\
& = R^{2} \int_{\Rbb} (1+|\xi|) |\widehat{(\chi\psi)}(R(\xi-\xi_0))|^2 d\xi \\
&= R\int_{\Rbb} \left( 1+ \left| \frac{\xi}{R} + \xi_0 \right | \right) |\widehat{(\chi\psi)}(\xi)|^2 d\xi \ge R \int_{\Rbb^1} |\widehat{(\chi\psi)}(\xi)|^2 d\xi.
\end{align*}
Thus we have \eqnref{fR1/2}.

On the other hand,
$$
\Fcal \left( \Gl f_R -\frac{1}{2}P_2 *f_R \right)(\xi)= \left( \Gl - \frac{1}{2} e^{-4\pi |\xi|} \right) R \widehat{(\chi\psi)}(R(\xi-\xi_0)).
$$
Thus,
\begin{align*}
\left\| \Gl f_R -\frac{1}{2}P_2 *f_R \right\|_{1/2}^2 &= R^{2} \int_{\Rbb^1} (1+|\xi|) \left| (\Gl - \frac{1}{2} e^{-4\pi |\xi|}) \widehat{(\chi\psi)}(R(\xi-\xi_0)) \right|^2 d\xi \\
&= R \int_{\Rbb^1} \left( 1+ \left| \frac{\xi}{R} + \xi_0 \right| \right) \left| \Gl - \frac{1}{2} e^{-4\pi | \frac{\xi}{R} + \xi_0 |} \right|^2 \left| \widehat{(\chi\psi)}(\xi) \right|^2 d\xi.
\end{align*}

Let
$$
\left\| \Gl f_R -\frac{1}{2}P_2 *f_R \right\|_{1/2}^2 = R \left( \int_{|\xi| \le \sqrt{R}} + \int_{|\xi| > \sqrt{R}} \right) =: I + II.
$$
If $|\xi| \le \sqrt{R}$, we have from \eqnref{GLdef}
$$
\left| \Gl - \frac{1}{2} e^{-4\pi | \frac{\xi}{R} + \xi_0 |} \right| \lesssim \frac{|\xi|}{R},
$$
and hence
\begin{align*}
I &= R \int_{|\xi| \le \sqrt{R}} \left( 1+ \left| \frac{\xi}{R} + \xi_0 \right| \right) \left| \Gl - \frac{1}{2} e^{-4\pi | \frac{\xi}{R} + \xi_0 |} \right|^2 \left| \widehat{(\chi\psi)}(\xi) \right|^2 d\xi \\
& \lesssim R^{-1} \int_{\Rbb^1} |\xi|^2 \left| \widehat{(\chi\psi)}(\xi) \right|^2 d\xi \lesssim R^{-1}.
\end{align*}

To estimate $II$, we observe that since $\chi\psi$ is a compactly supported smooth function,
$$
\left| \widehat{(\chi\psi)}(\xi) \right| \lesssim (1+|\xi|)^{-N}
$$
for any $N$. Thus, we have
$$
II \lesssim R \int_{|\xi| > \sqrt{R}} ( 1+ |\xi| )^{1-2N} d\xi \lesssim R^{1-N},
$$
we arrive at \eqnref{Gl-Kest}. \qed

\section{Proof of Theorem \ref{thm:main}}\label{sec3}

We prove the following proposition.

\begin{prop}\label{prop}
Let $\Gl \in (0, 1/2]$. There is a sequence $\Gvf_R \in H^{1/2}(\p\GO_{R})$ such that
\beq\label{limit}
\lim_{R \to \infty} \frac{\| ( \Gl I - \Kcal_{R}) [\Gvf_R] \|_{H^{1/2}(\p\GO_{R})}}{\| \Gvf_R \|_{H^{1/2}(\p\GO_{R})}} =0.
\eeq
\end{prop}

Theorem \ref{thm:main} is an immediate consequence of Proposition \ref{prop}. In fact, if $\Gl \in (0, 1/2]$, but $\Gl \notin \ol{\cup_{j=1}^\infty \Gs(\Kcal_{R_j})}$, then
$$
\mbox{dist} \left( \Gl, \ol{\cup_{j=1}^\infty \Gs(\Kcal_{R_j})} \right) >0.
$$
Thus, there is a constant $C$ independent of $j$ such that
$$
\| \Gvf \|_{H^{1/2}(\p\GO_{R_j})} \le C \| ( \Gl I - \Kcal_{R_j})[\Gvf] \|_{H^{1/2}(\p\GO_{R_{j}})}
$$
for all $\Gvf \in H^{1/2}(\p\GO_{R_{j}})$. Therefore, existence of the sequence $\{ \Gvf_R \}$ satisfying \eqnref{limit} implies that $\Gl \in \ol{\cup_{j=1}^\infty \Gs(\Kcal_{R_j})}$. Thus, $[0, 1/2] \subset \ol{\cup_{j=1}^\infty \Gs(\Kcal_{R_j})}$. Since the spectrum of the NP operator in two dimensions is symmetric with respect to $0$, we have $[-1/2, 1/2] \subset \ol{\cup_{j=1}^\infty \Gs(\Kcal_{R_j})}$. Since $\Gs(\Kcal_{R_j}) \subset [-1/2, 1/2]$, we arrive at \eqnref{main}.

\medskip
\noindent{\sl Proof of Proposition \ref{prop}}.
With the function $f_R$ defined in the previous section, define $\Gvf_R$ on $\p\GO_R$ by
\beq\label{Gvfdef3D}
\Gvf_R(x)= \Gvf_R(x_1,x_2) :=
\begin{cases}
f_R(x_1) \quad &\mbox{if } x \in \GG_R^+ \cup \GG_R^-, \\
0 \quad &\mbox{if } x \in \GG_R^s.
\end{cases}
\eeq
Then, \eqnref{fR1/2} yields
\beq\label{Gvflower}
R^{1/2} \lesssim \| \Gvf_R \|_{H^{1/2}(\p\GO_R)} .
\eeq
We show that
\beq\label{KRbounded}
\| (\Gl I- \Kcal_{R})[\Gvf_R] \|_{H^{1/2}(\p\GO_R)} \lesssim 1.
\eeq
These two estimates yield \eqnref{limit}.

To prove \eqnref{KRbounded}, choose a constant $C>0$ so that $\GG_R^l \subset \{ (x_1,x_2) \, : \, x_1 < -R+C \}$ and
$\GG_R^r \subset \{ (x_1,x_2) \, : \, x_1 > R-C \}$. Note that we can choose such a constant independently of $R$ if $R$ is sufficiently large.
Let $\Gz_1(x_1,x_2)=\Gz_1(x_1)$ be a smooth function supported in $(-R+C, R-C)$ such that $\Gz_1=1$ on $[-R+2C, R-2C]$ assuming that $R$ is sufficiently large, and let $\Gz_2:=1-\Gz_1$. Then we have
$$
\| (\Gl I- \Kcal_{R})[\Gvf_R] \|_{H^{1/2}(\p\GO_R)} \le \sum_{j=1}^2 \| \Gz_j(\Gl I- \Kcal_{R})[\Gvf_R] \|_{H^{1/2}(\p\GO_R)} .
$$

Thanks to \eqnref{Poisson}, we see that if $x \in \GG_R^+ \cup \GG_R^-$, then
\beq
\Kcal_{R} [\Gvf_R](x) = \frac{1}{2} (P_2 * f_R)(x_1).
\eeq
Thus, we have
$$
\Gz_1(\Gl I- \Kcal_{R})[\Gvf_R](x) = \Gz_1(x_1) \left( \Gl f_R(x_1) -\frac{1}{2} (P_2 * f_R)(x_1) \right).
$$
It then follows from \eqnref{Gl-Kest} that
\beq
\| \Gz_1(\Gl I- \Kcal_{R})[\Gvf_R] \|_{H^{1/2}(\p\GO_R)} \lesssim R^{-1/2}.
\eeq

We now estimate $\| \Gz_2(\Gl I- \Kcal_{R})[\Gvf_R] \|_{H^{1/2}(\p\GO_R)}$. Let
$$
\GG:= \p \GO_R \cap \{ (x_1,x_2) \, : \, x_1 < -R+C \mbox{ or } x_1 > R-C \}.
$$
Then the shape of $\GG$ is independent of $R$. Note that
$$
\Gz_2(\Gl I- \Kcal_{R})[\Gvf_R] = \Gz_2 \Kcal_{R}[\Gvf_R]
$$
and
$$
\| \Gz_2(\Gl I- \Kcal_{R})[\Gvf_R] \|_{H^{1/2}(\p\GO_R)} = \| \Gz_2 \Kcal_{R}[\Gvf_R] \|_{H^{1/2}(\GG)}.
$$
We use the following characterization of the space $H^{1/2}(\GG)$ (see, e.g., \cite{GT}):
\beq
\| h \|_{H^{1/2}(\GG)}^2= \| h \|_{L^2(\GG)}^2 + \int_{\GG} \int_{\GG} \frac{|h(x)-h(z)|^2}{|x-z|^2} d\Gs(x) d\Gs(z).
\eeq

Let $k_R(x,y)$ be the integral kernel of $\Kcal_R$, namely,
$$
k_R(x,y) = \frac{1}{2\pi} \frac{\la y-x, {\bf \nu}_{y} \ra}{|x-y|^2}.
$$
If $x, z \in \mbox{supp}(\Gz_2)$ and $y \in \mbox{supp}(\Gvf_R)$, then
$$
|x-z| \lesssim 1, \quad |x-y| \gtrsim R, \quad |z-y| \gtrsim R.
$$
Thus, $|k_R(x,y)| \lesssim R^{-1}$. It thus follows from \eqnref{fRx} that
\beq\label{100}
\| \Gz_2 \Kcal_{\p\GO_R}[\Gvf_R] \|_{L^2(\GG)} \lesssim R^{-1} \int_{\Rbb^1} |f_R(x_1)|dx_1 \lesssim \int_{\Rbb^1} |(\chi\psi)(x_1)|dx_1 \lesssim 1.
\eeq
We also have
\begin{align*}
&|\Gz_2(x) k_R(x,y) - \Gz_2(z) k_R(z,y)| \\
& \le |\Gz_2(x)- \Gz_2(z)| |k_R(x,y)| + |\Gz_2(z)| |k_R(x,y)-k_R(z,y)| \\
& \lesssim R^{-1} |x- z| + R^{-2} |x- z| \le R^{-1} |x- z|.
\end{align*}
Thus we have
\begin{align*}
|\Gz_2(x) \Kcal_R[\Gvf_R](x) - \Gz_2(z) \Kcal_R[\Gvf_R](z)| &\lesssim R^{-1}|x- z| \int_{\Rbb^1} |f_R(x_1)|dx_1 \\
&\le |x- z| \int_{\Rbb^1} |(\chi\psi)(x')|dx' \\
& \lesssim |x- z|.
\end{align*}
It then follows that
$$
\int_{\GG} \int_{\GG} \frac{|\Gz_2(x) \Kcal_R[\Gvf_R](x) - \Gz_2(z) \Kcal_R[\Gvf_R](z)|^2}{|x- z|^2} ds(x) ds(z) \lesssim 1,
$$
which together with \eqnref{100} implies
$$
\| \Gz_2(\Gl I- \Kcal_R)[\Gvf_R] \|_{H^{1/2}(\p\GO_R)} \lesssim 1.
$$
Thus, \eqnref{KRbounded} follows. This completes the proof. \qed

\section*{Discussions}

The spectral property \eqnref{main} on thin domains of rectangular shape is shared by NP operators on thin ellipses. In fact, if $E_j$, $j=1,2, \ldots$, is the ellipse defined by $x_1^2/a_j^2 + x_2^2/b_j^2 <1$ and $\Kcal_{\p E_j}$ is the corresponding NP operator, where $a_j$ and $b_j$ are positive numbers such that $b_j<a_j$ for all $j$ and $b_j/a_j \to 0$ as $j \to \infty$, then
\beq\label{main2}
\ol{\cup_{j=1}^\infty \Gs(\Kcal_{\p E_j})} = [-1/2, 1/2].
\eeq
We present a short proof of this fact below for readers' sake. It would be quite interesting to characterize the geometric properties of the family of thin domains which guarantee the spectral properties like \eqnref{main} and \eqnref{main2}.

There are at least two different kinds of thin domains in three dimensions: thin plate-like domains and thin cylinder-like domains. In the first case, it can be proved by modifying the proof of this paper that the same spectral property of the NP operators holds. However, it seems that the NP operators on the second kind of thin domains exhibit a completely different spectral structure. This investigation is in progress and the outcome will be reported in a forthcoming paper.

Let us now prove \eqnref{main2}.
It is known that the spectrum $\Gs(\Kcal_{\p E_j})$ of the NP operator on $E_j$ is $\{ \pm \frac{1}{2} r_j^n \, : \, n=1,2,\ldots \}$, where $r_j = \frac{a_j-b_j}{a_j+b_j}$ (see \cite{Ahl}).
Thus,
$$
\cup_{j=1}^\infty \Gs(\Kcal_{\p E_j}) = \{ \pm \frac{1}{2} r_j^n \, : \, j=1,2, \ldots, \ n=1,2,\ldots \}.
$$
Let $t_j: = - \log r_j$, and consider the set $S:=\{ n t_j \, : \, j=1,2, \ldots, \ n=1,2,\ldots \}$. Note that $t_j>0$ and $t_j \to 0$ as $j \to \infty$ since $r_j \to 1$. One can easily show that for any pair of real numbers $0 \le p <q$, there are $n$ and $j$ such that $p < nt_j <q$. Thus, $S$ is dense in $[0, \infty)$. We then infer that $\{ r_j^n \, : \, j=1,2, \ldots, \ n=1,2,\ldots \}$ is dense in $[0,1]$. So, \eqnref{main2} follows.



\begin{thebibliography}{11}

\bibitem{Ahl} {L.V. Ahlfors}:  Remarks on the Neumann-Poincar\'e
integral equation, Pacific J. Math. 3 (1952),  271--280.

\bibitem{HKL} J. Helsing, H. Kang and M. Lim, Classification of spectra of the Neumann--Poincar\'{e} operator on planar domains with corners by resonance, Ann. I. H. Poincare-AN 34 (2017), 991--1011.

\bibitem{GT} D. Gilbarg and N. Trudinger, {\sl Elliptic Partial Differential Equations of Second Order}, Springer-Verlag Berlin Heidelberg, 2001.

\bibitem{KPS-ARMA-07} {D. Khavinson, M. Putinar,  and H. S. Shapiro},
    {Poincar\'e's variational problem in potential theory},
  {Arch. Ration. Mech. An.} 185 (2007), 143--184.

\bibitem{PP2} K.M. Perfekt and M. Putinar, The essential spectrum of the Neumann-Poincare operator on a domain with corners, Arch. Ration. Mech. An. 223 (2017), 1019-1033.

\bibitem{Yo} K. Yosida, {\sl Functional Analysis}, 4th Ed., Springer, Berlin, 1974.

\end{thebibliography}
\end{document}